\documentclass[11pt]{amsart}
\usepackage{amssymb,latexsym,comment,url,amsmath}

\usepackage[T1]{fontenc}
\usepackage[misc]{ifsym}
\usepackage{xcolor}

\newcommand{\Char}{\operatorname{char}}

\newcommand{\disc}{\operatorname{disc}}

\newcommand{\GL}{\operatorname{GL}}

\newcommand{\Q}{{\mathbb Q}}
\newcommand{\OO}{{\mathcal O}}

\newcommand{\Z}{{\mathbb Z}}

\newcommand{\Magma}{{\sf MAGMA}}
\newcommand{\Mathematica}{{\sf Mathematica}}

\newenvironment{Proof}{\par\noindent{\sc Proof:}}%
                      {\hspace*{\fill}\nobreak$\Box$\par\medskip}
                       {\hspace*{\fill}\nobreak$\Box$\par\medskip}

\newtheorem{Proposition}{Proposition}[section]
\newtheorem{Theorem}[Proposition]{Theorem}
\newtheorem{Lemma}[Proposition]{Lemma}
\newtheorem{Corollary}[Proposition]{Corollary}

\theoremstyle{definition}

\newtheorem{Remark}[Proposition]{Remark}
 
\renewcommand{\baselinestretch}{1.1}

\renewcommand{\thefootnote}{\fnsymbol{footnote}}

\begin{document}

\title[Genus $2$ curves with bad reduction at one odd prime]%
{Genus $2$ curves with bad reduction at one odd prime}

\author[A. D\k{a}browski]%
{Andrzej~D\k{a}browski}
\address{Institute of Mathematics, University of Szczecin, Wielkopolska 15, 70-451 Szczecin, Poland}
\email{dabrowskiandrzej7@gmail.com}
\email{andrzej.dabrowski@usz.edu.pl}

\author[M. Sadek]%
{Mohammad~Sadek$^{\dagger}$}
\address{Faculty of Engineering and Natural Sciences, Sabanc{\i} University, Tuzla, \.{I}stanbul, 34956 Turkey}
\email{mmsadek@sabanciuniv.edu}

\begin{abstract}{ The problem of classifying elliptic curves over $\mathbb Q$ with a given discriminant has received much attention. The analogous problem for genus $2$ curves has only been tackled  when the absolute discriminant is a power of $2$. In this article we classify genus $2$ curves $C$ defined over $\Q$ with at least two rational Weierstrass points and whose absolute discriminant is an odd prime. In fact, we show that such a curve $C$ must be isomorphic to a specialization of one of finitely many $1$-parameter families of genus $2$ curves.  In particular,  we provide genus $2$ analogues to Neumann-Setzer families of elliptic curves over the rationals.}
\end{abstract}

\maketitle

\let\thefootnote\relax\footnotetext{ 
\hskip-16pt$^{\dagger}$Corresponding author
The datasets generated during and/or analysed during the current study are available from the corresponding author on reasonable request.}

\section{Introduction}
A well known result of Shafarevich \cite{Sha} states that the number of isomorphism classes of elliptic
curves over a given number field that have good reduction outside a finite set of primes
is finite. The online tables by Cremona \cite{C} exhibit all elliptic curves over the rational field of conductors
up to $500000$, together with additional arithmetic data such as the torsion subgroup and the Mordell-Weil rank.
In \cite{CL}, Cremona and Lingham give an explicit algorithm
to find all elliptic curves over a number field with good reduction outside a given
finite set of primes.

We remark that all results concerning explicit classifications of elliptic curves over $\mathbb Q$ with bad reduction outside a finite set of primes $S$ target the case when $S$
consists of at most two primes.
In what follows we give short overview of such known results.
Such elliptic curves were completely classified when $S=\{2\}$ by Ogg \cite{Og},
and when $S=\{3\}$ by Hadano \cite{Ha}.
Setzer \cite{Set} classified all elliptic curves with prime conductor and a rational point of order $2$.
Ivorra \cite{Iv} classified elliptic curves over $\mathbb Q$ of conductor $2^kp$, where $p$ is an odd prime,
with a rational point of order $2$.  Bennett, Vatsal and Yazdani \cite{BVY} classified all elliptic curves
over $\mathbb Q$ with a rational $3$-torsion point and good reduction outside the set
$\{3,p\}$, for a fixed prime $p$. Furthermore, Howe \cite{Ho}, Sadek \cite{Sa} and
D\k{a}browski-J\k{e}drzejak \cite{DJ} studied the classification of elliptic curves over $\mathbb Q$ with good
reduction outside two distinct primes and with a rational point of fixed order $\geq 4$.
In addition, Best and Matschke \cite{BM} presented a database of elliptic curves
with good reduction outside the first six primes.

Shafarevich conjectured \cite{Sha} that for each number field $K$, finite set of places $S$, and integer $g\ge 2$,
there are only finitely
many $K$-isomorphism classes of curves of genus $g$ over $K$ with good reduction
outside $S$. The proof was sketched by him in the hyperelliptic case; for details see the papers by Parshin and Oort,  \cite{Oort,Parshin}.
Merriman and Smart \cite{MS} determined all curves of genus $2$ with a rational Weierstrass point and with good reduction
away from $2$, up to an equivalence relation which is
coarser than the relation of isogeny between the associated Jacobian varieties.
Smart \cite{Sma} produced an explicit list of all genus $2$ curves with good reduction
away from $2$ by transforming the problem into the problem of solving some $S$-unit
equations. Rowan \cite{Row} adapted the latter method in order to produce examples of genus $2$ curves with
good reduction away from the prime $3$. Recently, infinitely many examples of genus $2$ curves over quadratic fields were presented when $S$ is empty, more precisely, the authors furnished examples of genus $2$ curves defined over the rational field that attain everywhere good reduction after a quadratic base change, \cite{DS}. Genus $2$ Curve Search Results from LMFDB \cite{Booker,C} give many
(probably not all) genus $2$ curves
with absolute discriminant up to $10^6$, together with additional arithmetic information. An expository paper by Poonen \cite{Poonen}
contains some potential relevant projects.

It can be seen that genus two curves with good reduction away from an odd prime have not been studied thoroughly in literature. In this article we are interested in genus $2$ curves $C$ with $\mathbb Q$-rational Weierstrass
points. 
 We attempt to extend the existing lists of genus two curves in \cite{MS,Sma}, to include curves with bad reduction at only one prime different from $2$. The aim of this paper is to find explicitly genus two curves with $\mathbb Q$-rational Weierstrass points and with odd prime absolute discriminant. We assemble lists of such genus two curves, analogous to existing lists of elliptic curves with bad reduction at only one odd prime.

In this work we consider genus two curves $C$ that can be described by globally minimal Weierstrass equations over $\Q$ of the form $y^2+Q(x)y=P(x)$, where $\deg Q(x)\le 2$ and $P(x)$ is monic of degree $5$. Moreover, we assume that these curves possess at least two $\Q$-rational Weierstrass points. This implies that they can be described by integral equations of the form $y^2=xf(x)$, where $f(x)$ is monic of degree $4$. Moreover, the latter equation may be assumed to be minimal at every prime except at $2$. It turns out that if $f(x)$ is reducible, then the absolute discriminant of $C$ can never be an odd prime, except when $f(x)=(x-b)g(x)$ and $g(x)$ is irreducible. We show that there are many (conjecturally, infinitely many) genus $2$ curves $C$
defined by $y^2=x(x-b)g(x)$ (with $g(x)$ irreducible) and such that the discriminant of $C$ is $\pm p$,
where $p$ is an odd prime. 
Let us give two families of such curves. In fact, we will prove in \S\ref{sec:thm4} that these are the only families of such curves.

\begin{itemize}
\item[(i)]  Let $f(t) = 256t^4 - 2064t^3 + 4192t^2 + 384t - 1051$.   The hyperelliptic curve $C_t$ defined by the
(non-minimal) equation
\[
y^2 = x(x+1)(x^3+64tx^2+64(t+4)x+256),\qquad t\in\Z,
\]
has discriminant $\pm p$ for some odd prime $p$ if and only if $f(t)=\pm p$. One can easily check that for $0 < t < 100$, $f(t)$ is a prime exactly when
$$t\in\{3,4,5,7,13,20,26,31,40,42,43,46,48,51, 55,82,83,90,98\},$$ and for such values of $t$, the discriminant $\Delta_{C_t}=f(t)$. For instance, one has
 $\Delta_{C_3}=2837$,    $\Delta_{C_4}=997$,
$\Delta_{C_5}=7669$,   $\Delta_{C_7}=113749$,   $\Delta_{C_{13}}=3489397$, and
$\Delta_{C_{20}}=26131429$.

\item[(ii)] Let $g(t) = 256t^4 + 768t^3 - 800t^2 - 2064t - 6343$.  The hyperelliptic curve $C_t$ given by the
(non-minimal) equation
$$
y^2 = x(x-4)(x^3+(4t+1)x^2-4(4t+5)x+64), \qquad t\in\Z,
$$
has discriminant $\pm p$ for some odd prime $p$ if and only if $g(t)=\pm p$. For $0 < t < 100$, $g(t)$ is a prime exactly when
$$t\in\{3,6,10,12,13,18,23,25,27,31,35,44,51, 58,74,80,82,93,95\},$$ and for such values of $t$,
$\Delta_{C_t}=g(t)$ is an odd prime, e.g., $\Delta_{C_3}=21737$,    $\Delta_{C_6}=450137$,
$\Delta_{C_{10}}=3221017$,   $\Delta_{C_{12}}=6489209$,   $\Delta_{C_{13}}=8830537$, and
$\Delta_{C_{18}}=31050137$.
\end{itemize}

 Conjecturally, each of the above $2$ families contains
infinitely many genus $2$ curves of prime discriminant. Such a statement follows from the above discussion, and a
classical conjecture by Bouniakovsky \cite{Bou} concerning prime values of irreducible
polynomials $f(x)\in\mathbb Z[x]\colon$ if the set of values $f(\mathbb Z^{+})$ has no common
divisor larger than $1$, then $|f(x)|$ represents infinitely many prime numbers.
It is not difficult to give examples with very large discriminants, for instance,
$f(49983)=\Delta_{C_{49983}}=1597567383051905525717$ and
$f(69945)=\Delta_{C_{69945}}=6126558731378331096629$ are primes,
where $f(t) = 256t^4 - 2064t^3 + 4192t^2 + 384t - 1051$, and $C_t$ belongs to the family (i) above.

In \S\ref{sec:thm5}, we give two explicit (conjecturally, infinite) families of genus $2$ curves with absolute prime discriminant described by $y^2=xf(x)$, with $f(x)$ an irreducible monic polynomial. We remark that the fact that we are looking for Weierstrass equations with odd prime absolute discriminant describing these curves implies that these Weierstrass equations are globally minimal.

It is worth mentioning that the families of genus $2$ curves that we obtain in this work can be seen as the genus-2 analogue to the famous Neumann-Setzer families of elliptic curves over the rationals \cite{Set}.
We recall that a Neumann-Setzer elliptic curve possesses a rational point of order $2$ and its  discriminant is an odd prime. Moreover, these elliptic curves may be described by the following globally minimal Weierstrass equation
$$y^2+xy=x^3+\frac{1}{4}(t-1)x^2-x,\qquad t\equiv 1 \textrm{ mod }4,$$
where the discriminant is an odd prime $p$ if and only if $t^2+64=p$, hence it is conjectured that there are infinitely many such curves.

Our explicit families of genus $2$ curves with odd prime (or odd square-free) discriminants lead to
abelian surfaces (Jacobians) with trivial endomorphisms, and may be useful when testing the paramodular conjecture
of Brumer and Kramer. If $C$ is such a curve, then the conjecture of Brumer and Kramer
predicts the existence of a cuspidal, nonlift Siegel paramodular
newform $f$ of degree $2$, weight $2$, and level $N_C$ with rational Hecke eigenvalues, such that
$L(\text{Jac}(C),s) = L(f,s,\text{spin})$. The interested reader may consult \cite{BK}.

 \section{Preliminaries on genus $2$ curves}

 Let $C$ be a smooth projective curve of genus $2$ over a perfect field $K$. Let $\sigma$ be the hyperelliptic involution of $C$. Given a generator $x$ of the subfield of $K(C)$ fixed by $\sigma$ over $K$, and $y\in K(C)$ such that $K(C)=K(x)[y]$, a {\em Weierstrass equation} $E$ of $C$ is given by
 \[E\colon y^2+Q(x)y=P(x),\qquad P(x),Q(x)\in K[x],\; \deg Q(x)\le 3,\;\deg P(x)\le 6.\]

 If $E'\colon v^2+Q'(u)v=P'(u)$ is another Weierstrass equation describing $C$, then there exist $\Big(
                                                                                           \begin{array}{cc}
                                                                                             a & b \\
                                                                                             c & d \\
                                                                                           \end{array}
                                                                                         \Big)\in\GL_2(K)
 $, $e\in K\setminus\{0\}$, $H(x)\in K[x]$ such that \[u=\frac{ax+b}{cx+d},\qquad v=\frac{ey+H(x)}{(cx+d)^3}.\]
If $\Char K\ne 2$, then we define the discriminant $\Delta_E$ of the Weierstrass equation $E$ to be
\[\Delta_E=2^{-12}\disc(4P(x)+Q(x)^2).\] One has $\Delta_E\ne 0$ if and only if $E$ describes a smooth curve. Moreover, \begin{eqnarray}\label{tr1}\Delta_{E'}=e^{20}(ad-bc)^{-30}\Delta_{E},\end{eqnarray} see for example \cite[\S 2]{Liu}.

 Assuming, moreover, that $K$ is a discrete valuation field with discrete valuation $\nu$ and ring of integers $\OO_K$, $E$ is said to be an {\em integral} Weierstrass equation of $C$ if both $P(x),Q(x)\in\OO_K[x]$. This implies that $\Delta_E\in\OO_K.$ A Weierstrass equation $E$ describing $C$ is said to be {\em minimal} if $E$ is integral and $\nu(\Delta_E)$ is the smallest valuation among all integral Weierstrass equations describing $C$. In the latter case, $\nu(\Delta_E)$ is the {\em discriminant} of $C$ over $\OO_K$.

 If $K$ is a number field with ring of integers $\OO_K$, then a Weierstrass model $E$ describing $C$ is {\em integral} if $P(x),Q(x)\in\OO_K[x].$ A Weierstrass equation $E$ is {\em globally minimal} if it is minimal over $\OO_{K_{\mathfrak p}}$ for every prime ideal $\mathfrak p$ of $\OO_K$, where $K_{\mathfrak p}$ is the completion of $K$ at $\mathfrak p$. Globally minimal Weierstrass equations do not exist in general, yet if $K$ has class number one, then $C$ has a globally minimal Weierstrass equation, \cite[Remarque 6]{Liu}. In the latter case, the discriminant of a globally minimal Weierstrass equation describing $C$ is the {\em discriminant} of $C$.

One notices that since we will be looking for Weierstrass equations with odd prime absolute discriminant, it follows that these equations are globally minimal, hence the corresponding discriminants are minimal.   

 \section{Rational Weierstrass points}

 In this section we assume that $C$ is a smooth projective genus $2$ curve defined over a number field $K$ of class number one. We assume moreover that $C$ possesses a $K$-rational Weierstrass point. It follows that $C$ can be described by a Weierstrass equation of the form
 \begin{eqnarray}\label{eq1}E\colon y^2+Q(x)y=P(x), \qquad\textrm{where }P(x),Q(x)\in K[x]\end{eqnarray}
 and $\deg Q(x)\le 2$, and $P(x)$ is monic of degree 5.

 Moreover, such an
equation is unique up to a change of coordinates of the form $x\mapsto u^2 x+r$, $y\mapsto u^5y+H(x)$ where $u\in
K\setminus\{0\}$, $r\in K$, and $H(x)\in K[x]$ is of degree at most $2$,
 see \cite[Proposition 1.2]{Loc}.

 Throughout this paper we will assume that $C$ is defined over $\Q$ by a globally minimal Weierstrass equation $E$ of the form in (\ref{eq1}). After the following transformation $x\mapsto x$ and $y\mapsto y+Q(x)/2$, then $C$ is described by $4y^2=4P(x)+Q(x)^2$. Now using the transformation $x\mapsto x/2^2$, $y\mapsto y/2^5$, an integral Weierstrass equation describing $C$ is $E'\colon y^2=G(x)$ where $G(x)\in\Z[x]$ is monic of degree $5$ and $\Delta_{E'}=2^{40}\Delta_E$.
\begin{Lemma}
\label{lem1}
Let $C$ be a smooth projective curve of genus $2$ defined over $\Q$ by a globally minimal Weierstrass equation of the form $y^2+Q(x)y=P(x)$, where $\deg Q(x)\le 2$ and $P(x)$ is monic of degree 5, with odd discriminant $\Delta$. Assume moreover that $C$ has at least two $\Q$-rational Weierstrass points. Then $C$ can be described by a Weierstrass equation of the form $E\colon y^2=xF(x)$, where $F(x)\in\Z[x]$ is a monic polynomial of degree $4$, and $\Delta_E=2^{40}\Delta$. In particular, $E$ is minimal over every $p$-adic ring $\Z_p$ except when $p=2$.
\end{Lemma}
\begin{Proof}
This follows from the argument above together with the fact that one of the rational Weierstrass points is sent to infinity, while the other point is sent to $(0,0)\in C(\Q)$ via a translation map. We notice that all the transformations used do not change minimality at odd primes.
\end{Proof}
 Let $C$ be a smooth projective curve of genus $2$ defined by a Weierstrass equation of the form $E\colon y^2=P(x)$, where $P(x)\in\Z[x]$ is of degree $5$ (not necessarily monic). The Igusa invariants $J_{2i}$, $1\le i\le 5$, associated to $E$ were defined in \cite[\S~4]{Igusa}. In fact these invariants can be defined for any Weierstrass equation describing $C$, see \cite{Liu3}. These invariants can be used to identify the reduction type of $C$ at a given prime $p$, see \cite{Igusa,Liu2}. For instance, the following result is \cite[Th\'{e}or\`{e}me 1]{Liu2}.

 \begin{Theorem}
 \label{theorem1}Let $C$ be a smooth projective curve of genus $2$ defined by the Weierstrass equation $y^2+Q(x)y=P(x)$ over $\Q$. Then $C$ has potential good reduction at the prime $p$ if and only if $J_{2i}^5/J_{10}^i\in \Z_p$, for every $1\le i\le 5$, where $\Z_p$ is the ring of $p$-adic integers.
 \end{Theorem}

 One remarks that if $C$ does not have potential good reduction at a prime $p$, then $C$ does not have good reduction at $p$.

 \section{Curves with six rational Weierstrass points}

 We assume that $C$ is a smooth projective curve of genus $2$ over $\Q$. If $C$ has six $\Q$-rational Weierstrass points, then $C$ may be described by a Weierstrass equation of the form $$E\colon y^2=x(x-b_1)(x-b_2)(x-b_3)(x-b_4),\qquad b_i\in\Z,\,i=1,2,3,4.$$

\begin{Theorem}
\label{thm1}
Let $C$ be a smooth projective curve of genus $2$ defined over $\Q$. Assume that $C$ has six $\Q$-rational Weierstrass points.
If $C$ is described by a globally minimal Weierstrass equation $E$ such that $|\Delta_E|$ is of the form $2^ap^b$, where $p$ is an odd prime, $a\ge0$, $b\ge 1$, then
$C$ is isomorphic to one of the following curves described by the following Weierstrass equations:  
\begin{eqnarray*}  E_0&\colon& y^2 = x (x-1) (x+1) (x-2) (x+2),\qquad \Delta_{E_0}=2^{18} \cdot 3^4,\\
E_1&\colon&   y^2 = x (x-3) (x+3) (x-6) (x+6),\qquad \Delta_{E_1}=2^{18} \cdot 3^{14}.
\end{eqnarray*}
\end{Theorem}
 \begin{Proof}
 The curve $C$ can be described by an integral Weierstrass equation of the form $E\colon y^2=x(x-b_1)(x-b_2)(x-b_3)(x-b_4)$, where $E$ is minimal at every odd prime. The discriminant $\Delta_{E}$ of $E$ is
described by
\[\Delta_{E}=2^8b_1^2 (b_1 - b_2)^2 b_2^2 (b_1 - b_3)^2 (b_2 - b_3)^2 b_3^2 (b_1 - b_4)^2 (b_2 - b_4)^2 (b_3 -
   b_4)^2 b_4^2.\]
Now we assume that $\Delta_{E}=2^mp^{n}$ where $m\ge 8$, $n\ge 1$.

We claim that at least two of the $b_i$'s are even.  Assume on the contrary that
$b_1 = \pm p^{\alpha_1}$, $b_2 = \pm p^{\alpha_2}$, $b_3 = \pm p^{\alpha_3}$
($\alpha_1 \geq \alpha_2 \geq \alpha_3 \geq 0$) are all odd. Then
$|b_1 - b_2| = 2^{s_1}p^{l_1}$, $|b_1 - b_3| = 2^{s_2}p^{l_2}$, $|b_2 - b_3| = 2^{s_3}p^{l_3}$,
with $s_i \geq 1$, $i=1,2,3$.
If all $b_i$'s are of the same sign, then using Catalan's conjecture (Mih\u{a}ilescu's theorem) we obtain
$\alpha_1 = \alpha_2 + 2 = \alpha_3 +2$ and $\alpha_2 = \alpha_3 +2$, a contradiction.
Now, if some $b_i$ and $b_j$ are of opposite signs, then we obtain $\alpha_i = \alpha_j$, in particular, it follows that $\alpha_1 = \alpha_2 = \alpha_3$. This will imply that two of $b_i$'s are equal, which is a contradiction.

This justifies considering the following subcases:

(i) In case two of the $b_i$'s are even, we may assume without loss of generality that $b_1=\pm 2^{c_1}p^{d_1},b_2=\pm 2^{c_2}p^{d_2}$, $b_3=\pm p^{d_3}$,
$b_4=\pm p^{d_4}$. with $c_1\ge c_2>0$. Elementary, but long case by case calculations show that
necessarily we have $d_1=d_2=d_3=d_4=d$; in particular $b_3=-b_4$.
Now, it is easy to check, that  $p=3$ and $c_1=c_2=1$; in particular $b_1=-b_2$.
Hence $b_1=2\cdot 3^d$, $b_2=-2 \cdot 3^d$, $b_3=3^d$, $b_4=-3^d$, which leads to the
Weierstrass equation $E_d\colon   y^2 = x (x-2 \cdot 3^d) (x+2 \cdot 3^d) (x-3^d) (x+3^d)$. A straight forward change of variables yield that the Weierstrass equations
$E_d$ and $E_{d+2}$ describe two isomorphic genus $2$ hyperelliptic curves, hence we only obtain two non-isomorphic genus $2$ curves $C_0$ and $C_1$ in the latter family described by $E_0$ and $E_1$ with
 minimal discriminants $2^{18} \cdot 3^4$ and $2^{18} \cdot 3^{14}$ respectively.

(ii)   We assume now without loss of generality that $b_1=\pm 2^{c_1}p^{d_1},b_2=\pm 2^{c_2}p^{d_2}$, $b_3=\pm 2^{c_3} p^{d_3}$,
$b_4=\pm p^{d_4}$. with $c_1\ge c_2\geq c_3 >0$.  Again,  long case by case calculations show that
necessarily we have $d_1=d_2=d_3=d_4=d$.  In this case, we obtain $b_1=2^{3} \cdot 3^d$,
$b_2=-2^{2} \cdot 3^d$, $b_3=2 \cdot 3^d$, $b_4=-3^d$, which leads to the curves $C'_d$ described by the Weierstrass equations
$  y^2 = x (x-2^{3} \cdot 3^d) (x+2^{2} \cdot 3^d) (x-2 \cdot 3^d) (x+ 3^d)$.
Again, it is easily seen that the curves $C'_{d}$ and $C'_{d+2}$ are isomorphic. Moreover, using the function \texttt{IsIsomorphic(,)} in \Magma, we can verify that the curves $C_0$ and $C'_0$ are isomorphic, and that
the curves $C_1$ and $ C'_1$ are isomorphic.

(iii) We assume now without loss of generality that $b_1=\pm 2^{c_1}p^{d_1},b_2=\pm 2^{c_2}p^{d_2}$, $b_3=\pm 2^{c_3} p^{d_3}$,
$b_4=\pm 2^{c_4} p^{d_4}$, with $c_1\ge c_2\geq c_3 \geq c_4 >0$. Again, long case by case calculations show that
necessarily we have $d_1=d_2=d_3=d_4=d$. In this case, we obtain $b_1=2^{t+3} \cdot 3^d$,
$b_2=-2^{t+2} \cdot 3^d$, $b_3=2^{t+1} \cdot 3^d$, $b_4=-2^t \cdot 3^d$, which lead to the curves
$C_{t,d}$ described by $y^2 = x (x-2^{t+3} \cdot 3^d) (x+2^{t+2} \cdot 3^d) (x-2^{t+1} \cdot 3^d) (x+ 2^t \cdot 3^d)$.
Now the curves $C_{t,d}$ and $C_{t,d+2}$ are isomorphic. Moreover, we may check using the function \texttt{IsIsomorphic(,)} that the curves $C_{t,d}$ and $C_{t+1,d}$
are isomorphic. Therefore, we obtain only two non-isomorphic curves $C_{1,0}$ and $C_{1,1}$. Finally, in a similar fashion one notices that 
 $C_0$ and $C_{1,0}$ are isomorphic, and the genus $2$ curves $C_1$ and $ C_{1,1}$ are isomorphic.
\end{Proof}

 \begin{Remark}
 One sees easily that none of the curves $C$ described in Theorem \ref{thm1} can be described by a globally minimal Weierstrass equation whose discriminant is square-free. This holds because $\Delta_{E}$ is always a square. Moreover, if $C$ is a curve that is described by neither $E_0$ nor $E_1$, and $C$ has bad reduction at exactly two primes, then both primes must be odd.
 \end{Remark}

 \begin{Corollary}
 Let $C$ be a smooth projective curve of genus $2$ defined over $\Q$. Assume that $C$ has six $\Q$-rational Weierstrass points. If $C$ is described by a globally minimal Weierstrass equation $E$, then $|\Delta_E|$ can never be a power of a prime. In other words, $C$ cannot have bad reduction at exactly one prime.
\end{Corollary}

 \begin{Proof}
 Theorem \ref{thm1} asserts that if $C$ has bad reduction at exactly one prime, then this prime must be $2$. However, according to \cite[\S 6.1]{MS}, there is no such curve with bad reduction only at $2$.
 \end{Proof}

 \section{Curves with exactly four rational Weierstrass points}

 We assume that $C$ is a smooth projective curve of genus $2$ over $\Q$ described by a globally minimal Weierstrass equation of the form $E\colon  y^2+Q(x)y=P(x)$, $P(x),Q(x)\in\Z[x]$, $\deg Q(x)\le 2$, and $P(x)$ is monic of degree $5$. If $C$ has exactly four $\Q$-rational Weierstrass points, then $C$ may be described by a Weierstrass equation of the form $$E'\colon y^2 = x(x-b_1)(x-b_2)(x^2+b_3x+b_4),\quad b_i\in\Z,\;i=1,2,3,4,$$ with $\Delta_{E'}=2^{40}\Delta_E$, see Lemma \ref{lem1}.

 \begin{Theorem}
 \label{thm2}
 Let $C$ be a smooth projective curve of genus $2$ defined over $\Q$. Assume that $C$ has exactly four $\Q$-rational Weierstrass points. If $C$ is described by a globally minimal Weierstrass equation of the form $E\colon y^2+Q(x)y=P(x)$, $\deg Q(x)\le 2$ and $ P(x)$ is monic of degree $5$, then $|\Delta_E|$ is never an odd prime.
 \end{Theorem}

 \begin{Proof}
 In accordance with Lemma \ref{lem1}, $C$ is described by $E'\colon y^2 = x(x-b_1)(x-b_2)(x^2+b_3x+b_4)$, $b_i\in\Z$ and $x^2+b_3x+b_4$ is irreducible. Moreover, $\Delta_{E'}=2^{40}\Delta_E$, hence $E'$ is minimal at every odd prime.
We have the following explicit formula for the discriminant of $E'$:
\begin{equation}
\Delta_{E'} = 2^8 b_1^2(b_1-b_2)^2b_2^2(b_3^2-4b_4)b_4^2(b_1^2+b_1b_3+b_4)^2
(b_2^2+b_2b_3+b_4)^2.
\end{equation}
We now assume that $\Delta_{E'}=\pm 2^{40}p$, where $p$ is an odd prime. It follows that
\begin{itemize}
\item[(a)] $b_1=\pm 2^a$,  $b_2=\pm 2^b$,  $b_1-b_2=\pm 2^c$,  $b_4=\pm 2^d$,
\item[(b)] $b_3^2-4b_4=\pm 2^e p$ (note that $b_3^2-4b_4$ is the only non-square factor, hence it's the only one that
can be divisible by $p$),
\item[(c)] $b_1^2+b_1b_3+b_4=\pm 2^f$,  $b_2^2+b_2b_3+b_4=\pm 2^g$,
\end{itemize}
where $a,b,c,d,e,f,g$ are non-negative integers such that $2a+2b+2c+2d+e+2f+2g=32$.
We will consider the following three cases.

\bigskip

(i) $a=b=0$. Then necessarily $b_1=-b_2$, $c=1$,
and combining the equations (c), we obtain
$b_4=\pm 2^{f-1} \pm 2^{g-1} -1$ and $b_3=\pm 2^{f-1} \pm 2^{g-1}$. The first
one gives $1\pm 2^d = \pm 2^{f-1} \pm 2^{g-1}$.

If $d\geq 1$, then $f=1$ (and, therefore $g=d+1$) or $g=1$ (and, therefore $f=d+1$).
In this case $b_3$ is odd, and hence $e=0$, and we obtain $4d+6=32$, which is impossible.

If $d=0$, then $\pm 2^{f-1} \pm 2^{g-1} = 1 \pm 2^d = 2$ or $0$.
In the first case, $f=g=1$ and $b_3=0$ or $\pm 2$, and there are no $p$ satisfying (b).
In the second case, $f=g\geq 1$, and $b_3=0$ or $\pm 2^f$. In the last case, (b) implies
$e=2$, $4f=28$, and hence $p=2^{12} \pm 1$, which is not a prime.

\bigskip

(ii) $a=0$, $b\geq 1$. Then necessarily $b=1$ and $c=0$.
We obtain a contradiction, considering carefully all possible tuples $(d,e,f,g)$
satisfying $2d+e+2f+2g=30$, and combining the equations (b) and (c).

\bigskip

(iii) $a, b\geq 1$. Then $a=b$, $b_1=-b_2$ and $c=a+1$.  We have
$2a+2b+2c = 6a+2$, hence we have five cases to consider: $a=b\leq 5$.
For each such $a$, we consider $d\geq 0$, and try to find $e$, $f$ and $g$
using (b) and (c).  None of these cases lead to genus $2$ curve $E$ with odd prime
value of $|\Delta_E|$. We omit the details.
 \end{Proof}

 \section{Curves with exactly two rational Weierstrass points and a quadratic Weierstrass point}

  Let $C$ be a smooth projective curve of genus $2$ over $\Q$ described by a globally minimal Weierstrass equation of the form $E\colon y^2+Q(x)y=P(x)$, $P(x),Q(x)\in\Z[x]$, $\deg Q(x)\le 2$, and $P(x)$ is monic of degree $5$. If $C$ has exactly two $\Q$-rational Weierstrass points and a quadratic Weierstrass point, then Lemma \ref{lem1} implies that $C$ is described by a Weierstrass equation of the form $$E'\colon y^2 = x(x^2+a_1x+a_2)(x^2+b_1x+b_2),\quad a_i, b_i\in\Z,$$ where both $x^2+a_1x+a_2$ and $x^2+b_1x+b_2$ are irreducible, and $\Delta_{E'}=2^{40}\Delta_E.$
 \begin{Theorem}
 \label{thm3}
 Let $C$ be a smooth projective curve of genus 2 defined over $\Q$. Assume that $C$ has exactly two $\Q$-rational Weierstrass points and a quadratic Weierstrass point. If $C$ is described by a globally minimal Weierstrass equation of the form $E\colon y^2+Q(x)y=P(x)$, $\deg Q(x)\le 2$ and $ P(x)$ is monic of degree $5$, then $|\Delta_E|$ is never an odd prime.
 \end{Theorem}
 \begin{Proof}
 As seen above, $C$ can be described by an integral Weierstrass equation $E'\colon y^2=x(x^2+a_1x+a_2)(x^2+b_1x+b_2)$ with $\Delta_{E'}=2^{40}\Delta_E$. In particular, $E'$ is minimal at every odd prime. We have the following explicit formula for the discriminant of $E$:
\begin{equation}
\Delta_{E'} = 2^8(a_1^2-4a_2)a_2^2(b_1^2-4b_2)b_2^2K^2,
\end{equation}
where $K=a_2^2-a_1a_2b_1+a_2b_1^2+a_1^2b_2-2a_2b_2-a_1b_1b_2+b_2^2$.
We assume that $\Delta_{E'}=\pm 2^{40}p$ where $p$ is an odd prime. It is clear that
$|a_2|=2^a$ and $|b_2|=2^b$, with $a,b\geq 0$. Therefore we can assume
without loss of generality that

\begin{equation} \label{Eq1}
|a_1^2-4a_2|=2^c,
\end{equation}
and
\begin{equation} \label{Eq2}
|b_1^2-4b_2|=2^dp,
\end{equation}
where $c,d\geq 0$.  Note that $K$ is necessarily a power of $2$. We will solve
systems of these equations, controlling the condition $2a+2b+c+d+2v_2(K)=32$, where $v_2$ is the $2$-valuation.
We will consider the following four cases, with many subcases.

\bigskip

(i) $a+2=c$, and both $a$ and $c$ are even.  Note that

$(a,c)\in\{(0,2), (2,4), (4,6), (6,8), (8,10), (10,12)\}$.

\bigskip

(ii) $a+2=c$, and both $a$ and $c$ are odd. Note that

$(a,c)\in\{(1,3), (3,5), (5,7), (7,9), (9,11)\}$.

\bigskip

(iii) $a+2>c$, then necessarily $c$ is even.
Using \eqref{Eq1} we obtain that $a+2-c = 1$ or $3$. This gives rise to the following $11$ pairs $(a,c)$:

(iiia) $(1,2)$, $(3,4)$, $(5,6)$, $(7,8)$, $(9,10)$,

(iiib)  $(1,0)$, $(3,2)$, $(5,4)$, $(7,6)$, $(9,8)$, $(11,10)$.

\bigskip

(iv) $a+2<c$, then necessarily $a$ is even.
Using \eqref{Eq1} we obtain that $c-a-2 = 1$ or $3$. This yields the following $10$ pairs $(a,c)$:

(iva) $(0,3)$, $(2,5)$, $(4,7)$, $(6,9)$, $(8,11)$,

(ivb) $(0,5)$, $(2,7)$, $(4,9)$, $(6,11)$, $(8,13)$.

\bigskip

The general strategy of the proof is as follows:

\begin{itemize}

\item{} fix a pair $(a,c)$ as above (we have $32$ such pairs);

\item{}  we have $a_2=\pm 2^a$, hence we can calculate $a_1$ using \eqref{Eq1};

\item{} now consider $b_2=\pm 2^b$, for all non-negative integers $b$
satisfying $2a+2b+c\leq 32$. Then, of course, $d+2v_2(K) \leq 32-2a-2b-c$;

\item{} for each triple $(a,b,c)$ check whether there exist $b_1$ and $K$
satisfying \eqref{Eq2} and $d+2v_2(K) = 32-2a-2b-c$.  Here we use a more
 convenient expression for $K$, namely
$K=(a_2-b_2)^2+(a_1-b_1)(a_1b_2-a_2b_1)$.

\end{itemize}

\bigskip

The cases with large $2a+c$ are the easiest to consider, and the cases with small
$2a+c$ are the longest ones (many subcases, etc.). Let us illustrate the method
in one of the easiest cases, $(a,c)=(11,10)$.
Here we have $a_2=\pm 2^{11}$ and $a_1=\pm 2^53$. If $b_2=\pm 1$,
then $K=(2^{11}\pm 1)^2+(\pm 2^53-b_1)(\pm 2^53\pm 2^{11}b_1)
=   (2^{11}\pm 1)^2+2^5(\pm 2^53-b_1)(\pm 3\pm 2^{6}b_1) = \pm 1$.
Note that $b_1$ is odd (otherwise $d>0$ and $2a+2b+c+d>32$), hence the second
summand in $K$ is of the form $2^5s$, with odd $s$.
On the other hand, note  that $(2^{11}\pm 1)^2+1=2\times s_{\pm}$, and
 $(2^{11}\pm 1)^2-1=2^{12}\times t_{\pm}$, with odd $s_{\pm}$ and $t_{\pm}$,
a contradiction.  If $b_2$ is even, then $2a+2b+c>32$, again a contradiction.
 \end{Proof}

 We will discuss smooth curves of genus $2$ with exactly two rational Weierstrass points and no quadratic Weierstrass points separately.

 \section{Curves with exactly three rational Weierstrass points}\label{sec:thm4}

 In this section, we assume that $C$ is a smooth projective curve of genus $2$ over $\Q$ described by a globally minimal Weierstrass equation of the form $E \colon y^2+Q(x)y=P(x)$, where $P(x),Q(x)\in\Z[x]$, $\deg Q(x)\le 2$ and $\deg P(x)= 5$. Assume, moreover, that $\Delta_E$ is an odd square-free integer. In particular, $C$ has good reduction at the prime $2$. If, moreover, $C$ has exactly three $\Q$-rational Weierstrass points, then it follows from Lemma \ref{lem1} that $C$ is described by a Weierstrass equation of the form $$E'\colon y^2=x(x-b)(x^3+dx^2+ex+f),\qquad b,d,e,f\in\Z,$$ whose discriminant $\Delta_{E'}=2^{40}\Delta_E$, and such that $x^3+dx^2+ex+f$ is irreducible. This implies that $E'$ is minimal at every odd prime. In this section, we find explicitly all such genus $2$ curves. In fact, we show that there are only two one-parameter families of the latter Weierstrass equations.

 One has \begin{equation}
\Delta_{E'} = 2^8 b^2f^2(b^3+db^2+eb+f)^2 (d^2e^2-4e^3-4d^3f+18def-27f^2)=2^{40}\Delta_E,
\end{equation}
where $\Delta$ is an odd square-free integer.

Setting $\epsilon_i=\pm 1$, $i=1,2,3,4$, one has:
\begin{itemize}
\item[(a)] $b=\epsilon_1 2^k$, \; $f=\epsilon_2 2^l$,
\item[(b)] $b^3+db^2+eb+f = \epsilon_3 2^m$,
\item[(c)] $d^2e^2-4e^3-4d^3f+18def-27f^2 = \epsilon_4 2^n\Delta_E$,
\end{itemize}
where $2k+2l+2m+n=32$.

\begin{Theorem}
\label{thm4}
Let $C$ be a smooth projective curve of genus $2$ defined over $\Q$ with good reduction at the prime $2$. Assume that $C$ has exactly three $\Q$-rational Weierstrass points. If $C$ is described by a globally minimal Weierstrass equation of the form $E\colon y^2+Q(x)y=P(x)$, $\deg Q(x)\le 2$ and $ P(x)$ is monic of degree $5$, such that $|\Delta_E|$ is a square-free odd integer, then $E$ lies in one of the following two one-parameter globally minimal Weierstrass equations

\begin{itemize}

\item[(i)] $E_t\colon  y^2 - x^2\,y = x^5 + 16 t\, x^4 + (16 + 8 t )\,x^3 + (8 + t)\,x^2 + x$;

\item[(ii)] $F_t\colon y^2+(-x^2-x)\,y=x^5+(-1 + t)\,x^4 + (-2 - 2t)\,x^3+(2+ t)\,x^2 -  x$;

\end{itemize}

where $t\in\Z$.

\end{Theorem}

\begin{Proof}
As explained above, the curve $C$ can be described by an integral Weierstrass equation of the form $E'\colon y^2=x(x-b)(x^3+dx^2+ex+f)$, where $\Delta_{E'}=2^{40}\Delta_E$; and conditions (a), (b) and (c) are satisfied. The values of $b$ and $f$ are determined by (a). Condition (b) implies that $m\ge\min(k,l)$. If $l\ge k$, then
$e(t) = \epsilon_1\epsilon_3 2^{m-k}-\epsilon_1(2^{2k}\epsilon_1+2^{k}t+2^{l-k}\epsilon_2),$ where $d=t$. If $l<k$, then $m=l$; if $\epsilon_2=\epsilon_3$ then $e(t)=-(2^{2k}+2^k\epsilon_1 t)$ and $d=t$; whereas if $\epsilon_2=-\epsilon_3$ then $k=l+1$, $e(t)=-\epsilon_1\epsilon_2-\epsilon_1(2^{2k}\epsilon_1+2^kt)$ and $d=t$.
Therefore, in any case the Weierstrass equation $E_t'\colon=E'$ is described as follows \begin{equation}\label{eq2}E'_t \colon y^2=x(x-2^k\epsilon_1)(x^3+tx^2+e(t)x+2^l\epsilon_2),\quad t\in\Z.\end{equation}

The strategy of the proof now is as follows. Given a fixed pair of positive integers $(k,l)$ such that $0\le k+l\le 16$, $m$ is chosen such that $0\le m\le 16-k-l$, $m\ge \min(k,l)$, and $n=32-2k-2l-2m\ge 0$. One checks now which of these tuples $(k,l,m,n)$ yields a curve with good reduction at the prime $2$, given that condition (c) is satisfied, in particular \begin{equation}\label{eq3}2^n||(d^2e^2-4e^3-4d^3f+18def-27f^2).\end{equation} Let $E'_t$ be the corresponding integral Weierstrass equation, we first check whether it has potential good reduction at the prime 2. This can be accomplished using Theorem \ref{theorem1}. If it has potential good reduction at $2$, then one checks for which congruence classes of $t$, condition (\ref{eq3}) is satisfied.

In fact, the only Weierstrass equations $E'_t$ that describes a curve $C$ with potential good reduction at $2$, i.e., $J_{2i}^5/J_{10}^i\in \Z_2$, for every $1\le i\le 5$, and such that (\ref{eq3}) is satisfied are the ones corresponding to the following tuples $(k,l,m,n)$:

\begin{eqnarray*}(0,0,8,16),\, &\epsilon_1=-\epsilon_2,\quad t\equiv 3 \textrm{ mod }64,\\
(2,5,5,8), &t\equiv 2\textrm{ mod }4, \\
(1,6,3,12),\, &\epsilon_1=\epsilon_3, \quad t\equiv 0 \textrm{ mod }8,\\
(4,4,4,8),\, &\epsilon_2=\epsilon_3, \quad t\equiv 0 \textrm{ mod }4,\\
(2,6,6,4),\,  & t\equiv 1 \textrm{ mod }2,\\
(0,8,0,16),\, &\epsilon_1=\epsilon_3, \quad t\equiv 0 \textrm{ mod }64.
\end{eqnarray*}
Any other tuple $(k,l,m,n)$ will yield an integral Weierstrass equation for which $J_{2i}^5/J_{10}^i\not\in \Z_2$ for some $i$, $1\le i\le 5$; or condition (\ref{eq3}) is not satisfied by the corresponding Weierstrass equation. More precisely, any other tuple $(k,l,m,n)$ that is not in the above list yields an integral Weierstrass equation for which there is some $i$, $1\le i\le 5$, such that $J_{2i}^5/J_{10}^i=x_i(t)/y_i(t)$ where $x_i(t)-x_i(0)\in2\Z[t]$, $x_i(0)$ is an odd integer, and $y_i(t)\in 2\Z[t]$; or else it is impossible for $2^n$ to exactly divide $(d^2e^2-4e^3-4d^3f+18def-27f^2)$ for any choice of an integer value of $t$.

For the tuple $(2,6,6,4)$, the minimal discriminant equals $(16t^2+56t+157)^2$ if
$(\epsilon_1,\epsilon_2,\epsilon_3)=(1,1,-1)$, and it equals  $(16t^2-40t+133)^2$ if
$(\epsilon_1,\epsilon_2,\epsilon_3)=(-1,-1,1)$ (hence it is never square-free).

Note, that the models
$Y^2=X(X-\epsilon_1)(4X^3+(4t+2)X^2+2(-2\epsilon_1 t-2-\epsilon_1-\epsilon_1\epsilon_2+\epsilon_1\epsilon_3)X+2\epsilon_2)$ for $(2,5,5,8)$, and $Y^2=X(4X-2\epsilon_1)(X^3+2X^2+(-\epsilon_1 t-2\epsilon_1\epsilon_2)X+\epsilon_2)$  for
$(1,6,3,12)$ have  discriminants of the form $2^{20}\times \textrm{odd}$.
Such models are minimal at $2$, since the polynomials on the right hand side are twice a stable polynomial
(root multiplicities $< 3$) and it is not congruent to a square modulo $4$ (see \cite[Corollaire 2, p. 4594 ]{Liu} and
\cite{MSt}).

The tuple $(4,4,4,8)$, where $\epsilon_2=\epsilon_3=1$ and $t\equiv 0$ mod $4$, yields an integral Weierstrass equation $E_t'$ that defines a curve with good reduction at $2$ and $2^{40}||\Delta_{E_t'}$. Replacing $t$ with $4t$ and minimizing the equation $E_t'$ yields the curve described by
$$
E_t^1(\epsilon_1)\colon  y^2-x\,y=x^5+(-4 \epsilon_1 + t)\,x^4+(-16 - 8\epsilon_1 t)\,x^3+(64\epsilon_1+ 16 t)\,x^2-\epsilon_1\,x,
$$
with $2\nmid\Delta_{E_t^1}$ for any integer $t$.

The tuple $(0,8,0,16)$, where $\epsilon_1=\epsilon_3=-1$ and $t\equiv 0$ mod $64$, yields an integral Weierstrass equation $E_t'$ that defines a curve with good reduction at $2$ and $2^{40}||\Delta_{E_t'}$. Replacing $t$ with $64t$ and minimizing the equation $E_t'$ yields the equation
$$
E_t^2(\epsilon_2)\colon  y^2-x^2\,y=x^5+16 t\, x^4+(16 \epsilon_2 + 8 t )\,x^3+(8 \epsilon_2 + t)\,x^2+\epsilon_2\,x,
$$
with $2\nmid\Delta_{E_t^2}$ for any integer $t$.

The tuple $(0,0,8,16)$ where $\epsilon_1=1$ and $\epsilon_2=-1$, gives rise to an integral Weierstrass equation $E_t'$ that defines a curve with good reduction at $2$ and $2^{40}||\Delta_{E_t'}$, when $t\equiv 3$ mod $64$. Minimizing the equation $E_t'$ yields
$$
E_t^3(\epsilon_3)\colon   y^2+(-x^2-1)\,y=x^5+(-5+64\epsilon_3 - 16 t )\,x^4+ (9 - 208\epsilon_3 + 56 t )\,x^3+(-9 + 252 \epsilon_3- 73 t)\,x^2
$$
$$+
 (4 - 135\epsilon_3 + 42 t)\,x + (-1 + 27 \epsilon_3 - 9 t),
$$
such that $2\nmid\Delta_{E_t^3}$ for any integer $t$.

Now we can check, using \Magma, that the following tuples of Weierstrass equations describe isomorphic genus 2 curves:
$$
(E_t^2(1),\, E_{t-4}^2(-1),\, E_{t+4}^1(-1)); \quad (E_{t}^1(1),\, E_{-t}^1(-1)); \quad
(E_{t}^3(1),\, E_{-t}^2(-1)); \quad (E_{t}^3(-1),\, E_{-t}^2(1)).
$$

Similarly, for the tuple $(2,6,6,4)$ when $t\equiv 1$ mod $2$ and $(\epsilon_1,\epsilon_2,\epsilon_3)\not\in\{(1,1,-1),(-1,-1,1)\}$, this yields
$E_t^4(\epsilon_1,\epsilon_2,\epsilon_3)$:
$$
y^2+(-x^2-x)\,y=x^5+(-\epsilon_1 + t)\,x^4+(-3/2-\epsilon_1/2-\epsilon_1\epsilon_2 +
\epsilon_1\epsilon_3 - 2\epsilon_1 t)\,x^3+(\epsilon_1 +2 \epsilon_2 -\epsilon_3+ t)\,x^2-\epsilon_1\epsilon_2\,x
$$
after minimization where $2\nmid\Delta_{E_t^4(\epsilon_1,\epsilon_2,\epsilon_3)}$ for any integer $t$. Using \Magma, one checks that the following pairs of equations describe isomorphic genus 2 curves:

$$
E_t^4(1,-1,1) \textrm{ and } E_{t+1}^4(1,1,1); \quad
E_t^4(1,-1,-1) \textrm{ and } E_{t+2}^4(1,1,1); \quad
E_t^4(-1,1,1) \textrm{ and } E_{t-3}^4(1,1,1);
$$
$$
E_t^4(-1,1,-1) \textrm{ and } E_{t-2}^4(1,1,1); \quad
E_t^4(-1,-1,-1) \textrm{ and } E_{t-1}^4(1,1,1).
$$

Reasoning as in the cases of tuples $(2,5,5,8)$ and $(1,6,3,12)$, we obtain that,
in the remaining cases for the tuples $(0,0,8,16)$, $(4,4,4,8)$, $(2,6,6,4)$,  and $(0,8,0,16)$, the minimal discriminants are
of the form  $2^{20}\times \textrm{odd}$. 

After computing the discriminants for the above families of Weierstrass equations, one concludes that the Weierstrass equations for which the absolute value of the discriminant is a square-free odd integer lie only in the families $E_t:=E_t^2(1)$ and $F_t:=E_t^4(1,1,1)$.
\end{Proof}

\begin{Corollary}
The absolute discriminant $|\Delta_{E_{t_0}}|$  (resp. $|\Delta_{F_{t_0}}|$), $t_0\in\Z$, of the minimal Weierstrass equation
$E_{t_0}$ (resp. $F_{t_0}$)  is a square-free odd integer $m$ if and only if $|f(t_0)|=m$  (resp. $|g(t_0)|=m$)
where $f(t), g(t)\in\Z[t]$ are degree-$4$ irreducible polynomials described as follows

\begin{align*}
f(t) = 256 t^4  - 2064 t^3 + 4192 t^2 + 384 t  - 1051; \\
g(t) =   256 t^4 + 768 t^3 - 800 t^2 - 2064 t - 6343.
\end{align*}

In particular, $\Delta_{E_{t_0}}=\pm p$  (resp. $\Delta_{F_{t_0}}=\pm p$), $p$ is an odd prime,
if and only if $f(t_0)=\pm p$  (resp. $g(t_0)=\pm p$). It follows that there are,
conjecturally, infinitely many integer values $t$ such that $|\Delta_{E_{t}}|$ (resp. $|\Delta_{F_t}|$) is an odd prime. 
\end{Corollary}
\begin{Proof}
This follows immediately as direct calculations show that $\Delta_{E_t}=f(t)$ and $\Delta_{F_t}=g(t)$ where $E_t$ and $F_t$ are defined as in Theorem \ref{thm4}. Moreover, the polynomials $f(t)$ and $g(t)$ satisfy the conditions of Bouniakovsky' Conjecture, \cite{Bou}, for the infinitude of prime values attained by an irreducible polynomial.
\end{Proof}
Recall that $E_t^4(\epsilon_1,\epsilon_2,\epsilon_3)$ is the Weierstrass equation $$
y^2+(-x^2-x)\,y=x^5+(-\epsilon_1 + t)\,x^4+(-3/2-\epsilon_1/2-\epsilon_1\epsilon_2 +
\epsilon_1\epsilon_3 - 2\epsilon_1 t)\,x^3+(\epsilon_1 +2 \epsilon_2 -\epsilon_3+ t)\,x^2-\epsilon_1\epsilon_2\,x.
$$
The following statement is a corollary of the proof above.
\begin{Corollary}
Let $ (\epsilon_1,\epsilon_2,\epsilon_3)\in\{(1,1,-1),(-1,-1,1)\}$. There are, conjecturally, infinitely many integer values
$t$ such that $|\Delta_{E^4_t}(\epsilon_1,\epsilon_2,\epsilon_3)|=p^2$, $p$ is an odd prime.
\end{Corollary}

\section{Curves with exactly two rational Weierstrass points and no quadratic Weierstrass points}
\label{sec:thm5}

 Let $C$ be described by a Weierstrass equation of the form \[E\colon y^2=x (x^4 + bx^3 + cx^2 + dx+e),\qquad b,c,d,e\in\Z\]
where the quartic is irreducible. Then the discriminant is given by

 \begin{align*}
\Delta_E = &2^8 e^2 (b^2 c^2 d^2 - 4  c^3 d^2 - 4 b^3 d^3 + 18  b c d^3 -
   27  d^4 - 4 b^2 c^3 e \\
&+ 16  c^4 e + 18 b^3 c d e -
   80  b c^2 d e - 6  b^2 d^2 e + 144  c d^2 e \\
&- 27 b^4 e^2 +
   144  b^2 c e^2 - 128  c^2 e^2 - 192  b d e^2 + 256  e^3).
\end{align*}

In this section, although we were not able to utilize the methods used before to classify all such curves, we produce two one  parametric familes of curves that will contain infinitely many curves with an odd prime absolute discriminant. It is worth mentioning  that $\Delta_E/(2^4e^2)$ is the discriminant of the elliptic curve described by $E'\colon y^2=x^4+bx^3+cx^2+dx+e$, therefore  classifying genus $2$ curves with an odd prime absolute discriminant described by $E$ is equivalent to finding elliptic curves with odd  prime absolute discriminant defined by $E'$. 

\bigskip

(i)  Let $f(t)=6912t^4 + 6912t^3 + 2592t^2 + 432t - 65509$.  The hyperelliptic curve $C_t$ given by the
(non-minimal) equation
$$
y^2 = x(x^4+16(4t+1)x+256),\qquad t\in\Z,
$$
has discriminant $\pm p$ for some odd prime $p$ if and only if $f(t)=\pm p$. One can easily check that for $0 < t < 100$,
$|f(t)|$ is prime exactly when
$$t\in\{1,2,14,15,16,29,41,47,52,57,69,71,80,81\},$$ and for such values of $t$, the discriminant $\Delta_{C_t}=-f(t)$.
For instance, one has
 $\Delta_{C_2}=-111611$,
$\Delta_{C_{14}}=-284946491$,
$\Delta_{C_{15}}=-373772171$,
$\Delta_{C_{16}}=-481901339$,
$\Delta_{C_{29}}=-5059429931$, and
$\Delta_{C_{41}}=-20012351339$.
In a general case, $\Delta_{C_t}$ is an odd integer.

(ii)  Let $f(t)=6912t^4 - 19712t^3 + 167968t^2 - 288720t + 134075$.  The hyperelliptic curve $C_t$ given by the
(non-minimal) equation
$$
y^2 =x(x^4+(4t+1)x^3-80x^2+256x-256), \qquad t\in\Z,
$$
has discriminant $\pm p$ for some odd prime $p$ if and only if $f(t)=\pm p$. One can easily check that for $0 < t < 100$,
$f(t)$ is a prime exactly when
$$t\in\{1,4,7,14,36,39,44,67,81,96,99\},$$ and for such values of $t$, the discriminant $\Delta_{C_t}=-f(t)$.
For instance, one has
 $\Delta_{C_1}=-523$,
$\Delta_{C_{4}}=-2174587$,
$\Delta_{C_{7}}=-16177963$,
$\Delta_{C_{14}}=-240455387$,
$\Delta_{C_{36}}=-10897249403$, and
$\Delta_{C_{39}}=-15065561387$.
In a general case, $\Delta_{C_t}$ is an odd integer.

\bigskip

Conjecturally,  the above families contain
infinitely many genus $2$ curves with an odd prime absolute discriminant.

 \section*{Acknowledgement}
 The authors are very grateful to the anonymous referee for many corrections and valuable suggestions that improved the manuscript.
 The authors are very grateful to Armand Brumer and Ken Kramer, John Cremona, Qing Liu, and Michael Stoll
for useful correspondence and suggestions related to this work. All the calculations in this work were performed using \Magma \cite{Magma}, \Mathematica \cite{Mathematica}, and PARI/GP \cite{PARI2}. This work started while the authors were invited to Friendly Workshop on Diophantine Equations and Related Problems, 2019, at the Department of Mathematics, Bursa Uluda\v{g} University, Bursa-Turkey. The authors thank the colleagues of this institution for their hospitality and support. M. Sadek is partially supported by BAGEP Award of the Science Academy, Turkey, and The Scientific and Technological Research Council of Turkey, T\"{U}B\.{I}TAK; research grant: ARDEB 1001/122F312.

\end{document}